\documentclass[a4paper]{article}

\usepackage{amssymb,amsmath, amsfonts}
\usepackage{authblk}
\usepackage{pst-grad} 

\def\str{\rightarrow}

\def\mj{\mbox{\bf 1}}

\def\prop#1#2{\vspace{2ex} \noindent{\sc #1.} {\it #2} \par \vspace{2ex}}
\def\propo#1#2{\vspace{2ex} \noindent{\sc #1.} {\it #2}}
\def\dkz{\noindent{\sc Proof. }}
\def\df#1{\noindent{\sc Definition #1. }}
\def\qed{\hfill $\dashv$\vspace{2ex}}

\def\strt{\stackrel{\textbf{.}\,}{\rightarrow}}
\def\oW{\overline{W}}

\def\cM{{\cal M}}
\def\WM{\oW\!\cM}

\def\top{\mbox{\it{Top}}}
\def\cat{\mbox{\it{Cat}}}
\def\diag{\mbox{\rm{diag}}}

\begin{document}

\title{Segal's multisimplicial spaces}
\author{Zoran Petri\'c}
\affil{Mathematical Institute, SANU,\\ Knez Mihailova 36, p.f.\
367,\\ 11001 Belgrade, Serbia\\ \texttt{zpetric@mi.sanu.ac.rs} }

\date{}
\maketitle

\vspace{-3ex}

\begin{abstract}
\noindent Some sufficient conditions on a simplicial space
$X\!:\Delta^{op}\str\top$ guaranteeing that $X_1\simeq \Omega|X|$
were given by Segal. We give a generalization of this result for
multisimplicial spaces. This generalization is appropriate for the
reduced bar construction, providing an $n$-fold delooping of the
classifying space of a category.

\end{abstract}

\vspace{.3cm}


\vspace{.5ex}

\noindent {\small {\it Keywords$\,$}: simplicial space, loop
space, lax functor}

\vspace{.5ex}

\noindent {\small {\it Acknowledgements$\,$}: This work was
supported by project ON174026 of the Ministry of Education,
Science, and Technological Development of the Republic of Serbia.
}

\section{Introduction}
This note makes no great claim to originality. It provides a
complete inductive argument for a generalization of
\cite[Proposition~1.5]{S74}, which was spelled out, not in a
precise manner, in~\cite[paragraph preceding Theorem~2.1]{BFSV}.
The authors of \cite{BFSV} considered this generalization trivial
and did not even provide a sketch of a proof. Some related, but
quite different, results are given in \cite{C74} and \cite{B92}.

The main result of \cite{CP13} reaches its full potential role in
constructing a model for an $n$-fold delooping of the classifying
space of a category only with the help of such a generalization of
\cite[Proposition~1.5]{S74}. Although we referred to \cite{BFSV},
the referees of \cite{CP13} were not convinced that our bar
construction actually provides an appropriate model for delooping.
The aim of this note is to fill in a gap in the literature
concerning these matters.

Segal, \cite[Proposition~1.5]{S74}, gave conditions on a
simplicial space $X\!:\Delta^{op}\str\top$ guaranteeing that
$X_1\simeq \Omega|X|$. His intention was to cover a more general
class of simplicial spaces than we need for our purposes,
therefore he worked with nonstandard geometric realizations of
simplicial spaces. We generalize his result, in one direction, by
passing from simplicial spaces to multisimplicial spaces, but
staying in a class appropriate for the standard geometric
realization. Our result is formulated to be directly applicable to
the reduced bar construction of \cite{CP13}, providing an $n$-fold
delooping of the classifying space of a category.

We work in the category (here denoted by \top) of compactly
generated Hausdorff spaces. (This category is denoted by ${\cal
K}e$ in \cite{GZ67} and by \textbf{CGHaus} in \cite{ML71}.) The
objects of \top\ are called \emph{spaces} and the arrows are
called \emph{maps}. Products are given the compactly generated
topology. We adopt the following notation throughout: $\simeq$ for
homotopy of maps or same homotopy type of spaces, $\approx$ for
homeomorphism of spaces.

The category $\Delta$ (denoted by $\Delta^+$ in \cite{ML71}) is
the standard topologist's \emph{simplicial} category defined as in
\cite[VII.5]{ML71}. We identify this category with the subcategory
of \top\ whose objects are the standard ordered  simplices (one
for each dimension), i.e.,\ with the standard cosimplicial space
$\Delta\str\top$.

The objects of $\Delta$ are the nonempty ordinals $1,2,3,\ldots$,
which are rewritten as $[0],[1],[2],$ etc. The coface arrows from
$[n-1]$ to $[n]$ are denoted by $\delta_i^n$, for $0\leq i\leq n$,
and the codegeneracy arrows from $[n]$ to $[n-1]$ are denoted by
$\sigma^n_i$, for $0\leq i\leq n-1$.

For the opposite category $\Delta^{\!op}$, we denote the arrow
$(\delta_i^n)^{op}\!:[n]\str [n-1]$ by $d_i^n$ and
$(\sigma_i^n)^{op}\!:[n-1]\str [n]$ by $s_i^n$. For $f$ an arrow
of $\Delta^{op}$ (or $(\Delta^{op})^n$), we abbreviate $X(f)$ by
$f$ whenever the (multi)simplicial object $X$ is determined by the
context.

We consider all the monoidal structures to be strict, which is
supported by the strictification given by \cite[XI.3,
Theorem~1]{ML71}. Some proofs prepared for non-specialists are
given in the appendix.

\section{Multisimplicial spaces and their realization}

A \emph{multisimplicial space} is an object of the category
$\top^{(\Delta^{op})^n}$, i.e.,\ a functor from $(\Delta^{op})^n$
to \top. When $n=0$, this is just a space and when $n=1$, this is
a \emph{simplicial space}. As usual, for a multisimplicial space
$X:(\Delta^{op})^n\str \top$, we abbreviate
$X([k_1],\ldots,[k_n])$ by $X_{k_1\ldots k_n}$.

We say that $X:(\Delta^{op})^n\str \top$ is a
\emph{multisimplicial set} when every $X_{k_1\ldots k_n}$ is
discrete. A \emph{multisimplicial map} is an arrow of
$\top^{(\Delta^{op})^n}$, i.e.\ a natural transformation between
multisimplicial spaces. When $n=1$, this is a \emph{simplicial
map}. Throughout this paper we use the standard geometric
realization of (multi)simplicial spaces.

\vspace{2ex}

\df{2.1} The \emph{realization} functor
$|\mbox{\hspace{1ex}}|:\top^{\Delta^{op}}\str \top$ of simplicial
spaces is defined so that for a simplicial space $X$, we have
\[
|X|=\left(\coprod_n X_n\times\Delta^n\right){\Big\slash}\sim,
\]
where the equivalence relation $\sim$ is generated by
\[
(d^n_ix,t)\sim(x,\delta^n_it)\quad{\rm and}\quad
(s^n_ix,t)\sim(x,\sigma^n_it).
\]

\vspace{2ex}

\df{2.2} For $p\geq 0$, the functor
$\underline{\hspace{1.3ex}}^{(p)}:\mbox{\emph{Top}}^{(\Delta^{op})^{n+p}}\str
\mbox{\emph{Top}}^{(\Delta^{op})^n}$ of \emph{partial realization}
is defined inductively as follows

\vspace{.5ex}

$\underline{\hspace{1.3ex}}^{(0)}$ is the identity functor, and
$\underline{\hspace{1.3ex}}^{(p+1)}$ is the composition
\[
\mbox{\emph{Top}}^{(\Delta^{op})^{n+p+1}}\stackrel{\cong}{\longrightarrow}\left(\mbox{\emph{Top}}^{\Delta^{op}}\right)
^{(\Delta^{op})^{n+p}}\stackrel{|\mbox{\hspace{1ex}}|^{(\Delta^{op})^{n+p}}}{-\!\!\!-\!\!\!-\!\!\!-\!\!\!-\!\!\!-\!\!\!-\!\!\!\longrightarrow}
\mbox{\emph{Top}}^{(\Delta^{op})^{n+p}}\stackrel{\;\;\underline{\hspace{1.3ex}}^{(p)}}{\longrightarrow}\top^{(\Delta^{op})^n},
\]
where the first isomorphism maps $X$ to $Y$ such that
$(Y_{k_1\ldots k_{n+p}})_l=X_{k_1\ldots k_{n+p} l}$.

\vspace{2ex}

For a multisimplicial space $X:(\Delta^{op})^p\str \top$, we
denote $X^{(p)}$ by $|X|$. Hence, for $X:(\Delta^{op})^{n+p}\str
\top$, we have that $(X^{(p)})_{k_1\ldots k_n}=|X_{k_1\ldots
k_n\underbrace{\underline{\hspace{1.3ex}}\ldots
\underline{\hspace{1.3ex}}}_{p}}|$.

\vspace{2ex}

\df{2.3} If $X=Y^{(p)}$, for $Y$ a multisimplicial set, then we
say that $X$ is a \emph{partially realized} multisimplicial set
(PRmss).

\vspace{2ex}

\df{2.4} For $n\geq 0$ and $X:(\Delta^{op})^n\str \top$, let the
simplicial space $\diag X:\Delta^{op}\str \top$ be such that
\[
(\diag X)_k=X_{k\ldots k}.
\]
In particular, when $n=0$ and $X$ is just a topological space, we
have that $(\diag X)_k=X$ and all the face and degeneracy maps of
$\diag X$ are $\mj_X$.

\vspace{2ex}

The following lemma is a corollary of \cite[Lemma, p.\ 94]{Q73}.

\prop{Lemma 2.5}{For $X:(\Delta^{op})^n\str \top$, we have that
$|X|\approx|\diag X|$.}

As a consequence of Lemma~2.5 and \cite[Theorem~14.1]{M67} we have
the following lemma.

\prop{Lemma 2.6}{If $X$ is a PRmss, then $|X|$ is a CW-complex.}

The following remark easily follows.

\prop{Remark 2.7}{{\rm(a)} If $X:(\Delta^{op})^{n+p}\str \top$ is
a PRmss, then $X^{(p)}$ is a PRmss. \\[1ex] {\rm(b)} If $X:(\Delta^{op})^{n}\str \top$ is a
PRmss, then for every $k_1,\ldots,k_n$, the space $X_{k_1\ldots
k_n}$ is a CW-complex. \\[1ex] {\rm(c)} If $X:(\Delta^{op})^{n}\str \top$ is a
PRmss, then for every $k\geq 0$,
$X_{k\underline{\hspace{1.3ex}}\ldots \underline{\hspace{1.3ex}}}$
is
a PRmss. \\[1ex] {\rm(d)} If $X:(\Delta^{op})^{n}\str \top$, for $n>1$, is a
PRmss, then $Y\!:\Delta^{op}\times \Delta^{op}\str\top$, defined
so that $Y_{mk}=X_{mk\ldots k}$, is a PRmss}

\df{2.8} A simplicial space $X:\Delta^{op}\str \top$ is
\emph{good} when for every $0\leq i\leq n-1$, the map
$s^n_i\!:X_{n-1}\str X_n$ is a closed cofibration. It is
\emph{proper} (Reedy cofibrant) when for every $n\geq 1$, the
inclusion $sX_n\hookrightarrow X_n$, where $sX_n=\bigcup_i
s^n_i(X_{n-1})$, is a closed cofibration.

\prop{Proposition 2.9}{Every PRmss $X:\Delta^{op}\str \top$ is
good.}

\dkz Since $d^n_i\circ s^n_i=\mj_{X_{n-1}}$, we may consider
$X_{n-1}$ to be a retract of $X_n$. By Remark~2.7 (b), $X_n$ is a
CW-complex and by \cite[Corollary III.2]{DE72} (see also
\cite[Corollary 2.4 (a)]{L82}) a locally equiconnected space. By
\cite[Lemma 3.1]{L82} and since every $X_k$ is Hausdorff, $s^n_i$
is a closed cofibration. \qed

As a corollary of \cite[Proposition~22]{SR02} (see also references
listed in \cite[Section~6, p.\ 19]{SR02}) we have the following
result.

\prop{Lemma 2.10}{Every good simplicial space is proper.}

The following result is from \cite[Appendix, Theorem A4
(ii)]{M74}.

\prop{Lemma 2.11}{Let $f\!:X\str Y$ be a simplicial map of proper
simplicial spaces. If each $f_k\!:X_k\str Y_k$ is a homotopy
equivalence, then $|f|\!:|X|\str |Y|$ is a homotopy equivalence.}

\df{2.12} The \emph{product} $X\times Y$ of simplicial spaces $X$
and $Y$ is defined componentwise, i.e.\ $(X\times Y)_k=X_k\times
Y_k$, and for every arrow $f\!:k\str l$ of $\Delta^{op}$ and every
$x\in X_k$ and $y\in Y_k$, we have that $f(x,y)=(fx,fy)$.

\vspace{2ex}

Since the product of two CW-complexes in \top\ is a CW-complex, by
reasoning as in Proposition~2.9, we have the following.

\prop{Remark 2.13}{If $X,Y:\Delta^{op}\str \top$ are PRmss, then
$X\times Y$ is good.}

The following lemma is a corollary of \cite[Lemma~11.11]{M72}.

\prop{Lemma 2.14}{If the space $X_0$ of a simplicial space
$X\!:\Delta^{op}\str \top$ is path-connected, then $|X|$ is
path-connected.}

\section{Segal's multisimplicial spaces}

For $m\geq 1$, consider the arrows $i_1,\ldots, i_m:[m]\str [1]$
of $\Delta^{op}$ given by the following diagrams.

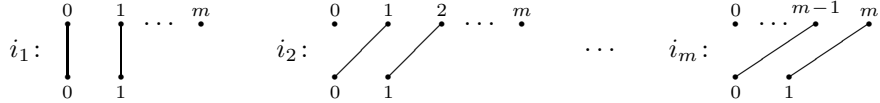
\begin{center}
\begin{picture}(320,40)(30,0)

\put(40,20){\makebox(0,0)[r]{$i_1\!:$}}

\put(50,10){\circle*{2}} \put(70,10){\circle*{2}}
\put(50,30){\circle*{2}} \put(70,30){\circle*{2}}
\put(85,30){\makebox(0,0){\ldots}} \put(100,30){\circle*{2}}

\put(50,2){\makebox(0,0)[b]{\scriptsize $0$}}
\put(70,2){\makebox(0,0)[b]{\scriptsize $1$}}
\put(50,33){\makebox(0,0)[b]{\scriptsize $0$}}
\put(70,33){\makebox(0,0)[b]{\scriptsize $1$}}
\put(100,33){\makebox(0,0)[b]{\scriptsize $m$}}

\put(50,10){\line(0,1){20}} \put(70,10){\line(0,1){20}}

\put(140,20){\makebox(0,0)[r]{$i_2\!:$}}

\put(150,10){\circle*{2}} \put(170,10){\circle*{2}}
\put(150,30){\circle*{2}} \put(170,30){\circle*{2}}
\put(190,30){\circle*{2}} \put(205,30){\makebox(0,0){\ldots}}
\put(220,30){\circle*{2}}

\put(150,2){\makebox(0,0)[b]{\scriptsize $0$}}
\put(170,2){\makebox(0,0)[b]{\scriptsize $1$}}
\put(150,33){\makebox(0,0)[b]{\scriptsize $0$}}
\put(170,33){\makebox(0,0)[b]{\scriptsize $1$}}
\put(190,33){\makebox(0,0)[b]{\scriptsize $2$}}
\put(220,33){\makebox(0,0)[b]{\scriptsize $m$}}

\put(150,10){\line(1,1){20}} \put(170,10){\line(1,1){20}}

\put(250,20){\makebox(0,0){\ldots}}

\put(290,20){\makebox(0,0)[r]{$i_m\!:$}}

\put(300,10){\circle*{2}} \put(320,10){\circle*{2}}
\put(300,30){\circle*{2}} \put(330,30){\circle*{2}}
\put(315,30){\makebox(0,0){\ldots}} \put(350,30){\circle*{2}}

\put(300,2){\makebox(0,0)[b]{\scriptsize $0$}}
\put(320,2){\makebox(0,0)[b]{\scriptsize $1$}}
\put(300,33){\makebox(0,0)[b]{\scriptsize $0$}}
\put(330,33){\makebox(0,0)[b]{\scriptsize $m\!-\!1$}}
\put(350,33){\makebox(0,0)[b]{\scriptsize $m$}}

\put(300,10){\line(3,2){30}} \put(320,10){\line(3,2){30}}

\end{picture}
\end{center}

The following images of these arrows under the functor ${\cal
J}\!:\Delta^{op}\str\Delta$ of \cite[Section~6]{PT13} may help the
reader to see that $i_1,\ldots,i_m$ correspond to $m$ projections.
(Note that 0 and 2 in the bottom line of the images serve to
project away all but one element of the top line.)

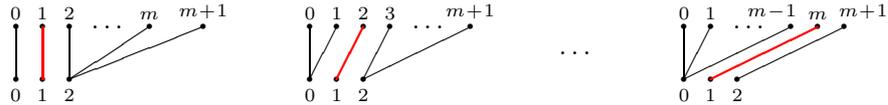
\begin{center}
\begin{picture}(320,40)(40,0)

\put(50,10){\circle*{2}} \put(60,10){\circle*{2}}
\put(70,10){\circle*{2}} \put(50,30){\circle*{2}}
\put(60,30){\circle*{2}} \put(70,30){\circle*{2}}
\put(85,30){\makebox(0,0){\ldots}} \put(100,30){\circle*{2}}
\put(120,30){\circle*{2}}

\put(50,2){\makebox(0,0)[b]{\scriptsize $0$}}
\put(60,2){\makebox(0,0)[b]{\scriptsize $1$}}
\put(70,2){\makebox(0,0)[b]{\scriptsize $2$}}
\put(50,33){\makebox(0,0)[b]{\scriptsize $0$}}
\put(60,33){\makebox(0,0)[b]{\scriptsize $1$}}
\put(70,33){\makebox(0,0)[b]{\scriptsize $2$}}
\put(100,33){\makebox(0,0)[b]{\scriptsize $m$}}
\put(120,33){\makebox(0,0)[b]{\scriptsize $m\!+\!1$}}

{\thinlines \put(50,10){\line(0,1){20}}
\put(70,10){\line(0,1){20}} \put(70,10){\line(3,2){30}}
\put(70,10){\line(5,2){50}}} {\thicklines
\put(60,10){\textcolor{red}{\line(0,1){20}}}}

\put(160,10){\circle*{2}} \put(170,10){\circle*{2}}
\put(180,10){\circle*{2}} \put(160,30){\circle*{2}}
\put(170,30){\circle*{2}} \put(180,30){\circle*{2}}
\put(190,30){\circle*{2}} \put(205,30){\makebox(0,0){\ldots}}
\put(220,30){\circle*{2}}

\put(160,2){\makebox(0,0)[b]{\scriptsize $0$}}
\put(170,2){\makebox(0,0)[b]{\scriptsize $1$}}
\put(180,2){\makebox(0,0)[b]{\scriptsize $2$}}
\put(160,33){\makebox(0,0)[b]{\scriptsize $0$}}
\put(170,33){\makebox(0,0)[b]{\scriptsize $1$}}
\put(180,33){\makebox(0,0)[b]{\scriptsize $2$}}
\put(190,33){\makebox(0,0)[b]{\scriptsize $3$}}
\put(220,33){\makebox(0,0)[b]{\scriptsize $m\!+\!1$}}

{\thinlines \put(160,10){\line(0,1){20}}
\put(160,10){\line(1,2){10}} \put(180,10){\line(1,2){10}}
\put(180,10){\line(2,1){40}}} {\thicklines
\put(170,10){\textcolor{red}{\line(1,2){10}}}}

\put(260,20){\makebox(0,0){\ldots}}

\put(300,10){\circle*{2}} \put(310,10){\circle*{2}}
\put(320,10){\circle*{2}} \put(300,30){\circle*{2}}
\put(310,30){\circle*{2}} \put(340,30){\circle*{2}}
\put(350,30){\circle*{2}} \put(325,30){\makebox(0,0){\ldots}}
\put(360,30){\circle*{2}}

\put(300,2){\makebox(0,0)[b]{\scriptsize $0$}}
\put(310,2){\makebox(0,0)[b]{\scriptsize $1$}}
\put(320,2){\makebox(0,0)[b]{\scriptsize $2$}}
\put(300,33){\makebox(0,0)[b]{\scriptsize $0$}}
\put(310,33){\makebox(0,0)[b]{\scriptsize $1$}}
\put(342,33){\makebox(0,0)[br]{\scriptsize $m\!-\!1$}}
\put(350,33){\makebox(0,0)[b]{\scriptsize $m$}}
\put(358,33){\makebox(0,0)[bl]{\scriptsize $m\!+\!1$}}

{\thinlines \put(300,10){\line(0,1){20}}
\put(300,10){\line(1,2){10}} \put(300,10){\line(2,1){40}}
\put(320,10){\line(2,1){40}}} {\thicklines
\put(310,10){\textcolor{red}{\line(2,1){40}}}}

\end{picture}
\end{center}

For maps $f_i\!:A\str B_i$, $1\leq i\leq m$, we denote by $\langle
f_1,\ldots,f_m\rangle:A\str B_1\times\ldots\times B_m$ the map
obtained by the Cartesian structure of \top. In particular, for
the above-mentioned $i_1,\ldots,i_m$ and for a simplicial space
$X\!:\Delta^{op}\str\top$ we have the map
\[
p_m=\langle i_1,\ldots,i_m\rangle\!:X_m\str (X_1)^m.
\]
(According to our convention from the introduction, $X(i_k)$ is
abbreviated by $i_k$.) If $m=0$, then $(X_1)^0=\{\ast\}$ (a
terminal object of \top) and let $p_0$ denote the unique arrow
from $X_0$ to $(X_1)^0$. The following lemma is claimed in
\cite{S74}.

\prop{Lemma 3.1}{If $X\!:\Delta^{op}\str\top$ is a simplicial
space such that for every $m\geq 0$, the map $p_m$ is a homotopy
equivalence, then $X_1$ is a homotopy associative
\mbox{\emph{H}-space} whose multiplication $m$ is given by the
composition
\[(X_1)^2\stackrel{p_2^{-1}}{\longrightarrow}X_2\stackrel{d^2_1}{\longrightarrow}X_1,\]
where $p_2^{-1}$ is an arbitrary homotopy inverse to $p_2$, and
whose unit $\ast$ is $s^1_0(x_0)$, for an arbitrary $x_0\in X_0$.}

\df{3.2} We say that a PRmss $X\!:\Delta^{op}\str\top$ is
\emph{Segal's simplicial space} when for every $m\geq 0$, the map
$p_m\!:X_m\str (X_1)^m$ is a homotopy equivalence.

\prop{Lemma 3.3}{Let $Y\!:\Delta^{op}\times \Delta^{op}\str\top$
be a PRmss. If for every $k\geq 0$, the simplicial space
$Y_{\underline{\hspace{1.3ex}}k}$ is Segal's, then $Y^{(1)}$ is
Segal's simplicial space.}

\df{3.4} We say that a PRmss $X\!:(\Delta^{op})^n\str\top$, where
$n\geq 1$, is \emph{Segal's multisimplicial space}, when for every
$l\in\{0,\ldots, n-1\}$ and every $k\geq 0$, the simplicial space
$X_{\underbrace{\mbox{\scriptsize $1\ldots
1$}}_{l}\underline{\hspace{1.3ex}}k\ldots k}$ is Segal's.

Note that we do not require $X_{k_1\ldots k_l
\underline{\hspace{1.3ex}} k_{l+1}\ldots k_{n-1}}$ to be Segal's
for arbitrary $k_1,\ldots,k_{n-1}$ (see the parenthetical remark
in Section~5.)

\propo{Remark 3.5}{If $X\!:(\Delta^{op})^n\str\top$ is Segal's
multisimplicial space, then for every $l\in\{0,\ldots, n-1\}$,
$X_{1\ldots 1}$ is homotopy associative \emph{H}-space with
respect to the structure obtained from Lemma~3.1 applied to
$X_{\underbrace{\mbox{\scriptsize $1\ldots
1$}}_{l}\underline{\hspace{1.3ex}} 1\ldots
1}\!:\Delta^{op}\str\top$.}

Our goal is to generalize the following proposition, which stems
from \cite[Proposition 1.5 (b)]{S74}. (In the proof of that
result, contractibility of $|PA|$ comes from the fact that
$|PA|\simeq A_0$.)

\prop{Proposition 3.6}{Let $X\!:\Delta^{op}\str\top$ be Segal's
simplicial space. If $X_1$ with respect to the \emph{H}-space
structure obtained by Lemma~3.1 is grouplike, then $X_1\simeq
\Omega|X|$.}

Our generalization is the following.

\prop{Proposition 3.7}{Let $X\!:(\Delta^{op})^n\str\top$ be
Segal's multisimplicial space. If $X_{1\ldots 1}$, with respect to
the \emph{H}-space structure obtained by Remark~3.5  when $l=n-1$,
is grouplike, then $X_{1\ldots 1}\simeq \Omega^n|X|$.}

\dkz We proceed by induction on $n\geq 1$. If $n=1$, the result
follows from Proposition~3.6.

If $n>1$, then we may apply the induction hypothesis to
$X_{1\underline{\hspace{1.3ex}}\ldots
\underline{\hspace{1.3ex}}}$. Hence,
\[
X_{1\ldots 1}\simeq
\Omega^{n-1}|X_{1\underline{\hspace{1.3ex}}\ldots
\underline{\hspace{1.3ex}}}|.
\]

By Lemma~2.5, we have that $|X_{1\underline{\hspace{1.3ex}}\ldots
\underline{\hspace{1.3ex}}}|\approx |\diag
X_{1\underline{\hspace{1.3ex}}\ldots
\underline{\hspace{1.3ex}}}|$. By the assumption and Remark
2.7~(d), the multisimplicial space $Y\!:\Delta^{op}\times
\Delta^{op}\str\top$, defined so that $Y_{mk}=X_{mk\ldots k}$,
satisfies the conditions of Lemma~3.3. Let $Z$ be the simplicial
space $Y^{(1)}\!:\Delta^{op}\str\top$, i.e.,\
\[
Z_m=|Y_{m\underline{\hspace{1.3ex}}}|= |\diag
X_{m\underline{\hspace{1.3ex}}\ldots \underline{\hspace{1.3ex}}}|.
\]

By Lemma 3.3, $Z$ is Segal's simplicial space. By Remark~2.7 (b),
$Z_1$ is a CW-complex. Since the space $Y_{10}$ (i.e.,
$X_{10\ldots 0}$) is by the assumption homotopic to $(X_{110\ldots
0})^0$, it is contractible, and hence, path-connected. By
Lemma~2.14, we have that $Z_1$, which is equal to
$|Y_{1\underline{\hspace{1.3ex}}}|$, is path-connected. Note also
that $|Z|=|Y|\approx |\diag X| \approx |X|$.

By Lemma~3.1, $Z_1$ is a homotopy associative H-space, and since
it is a path-connected CW-complex, by
\cite[Proposition~8.4.4]{A11}, it is grouplike. Applying
Proposition~3.6 to $Z$, we obtain
\[
|X_{1\underline{\hspace{1.3ex}}\ldots
\underline{\hspace{1.3ex}}}|\approx |\diag
X_{1\underline{\hspace{1.3ex}}\ldots
\underline{\hspace{1.3ex}}}|=Z_1 \simeq \Omega|Z| \approx
\Omega|X|.
\]
Finally, we have
\[
X_{1\ldots 1}\simeq
\Omega^{n-1}|X_{1\underline{\hspace{1.3ex}}\ldots
\underline{\hspace{1.3ex}}}| \simeq \Omega^n|X|.
\]

\vspace{-4ex}

\mbox{\hspace{1em}}\qed

\section{Segal's lax functors}

Thomason, \cite{T79}, was the first who noticed that the reduced
bar construction based on a symmetric monoidal category produces a
lax, instead of an ordinary, functor. The idea to use Street's
rectification in that case, also belongs to him.

We use the notions of \emph{lax functor}, \emph{left} and
\emph{right lax transformation} as defined in \cite{S72}. The
following theorem is taken over from \cite[Theorem~2]{S72}.

\prop{Theorem 4.1}{For every lax functor $W\!:{\cal A}\str \cat$
there exists a genuine functor $V\!:{\cal A}\str \cat$, a left lax
transformation $E\!:V\str W$ and a right lax transformation
$J\!:W\str V$ such that $J$ is the left adjoint to $E$ and
$W=EVJ$.}

\noindent We call $V$ a \emph{rectification} of $W$. It is easy to
see that if $W\!:{\cal A}\times{\cal B}\str \cat$ is a lax functor
and $V$ is its rectification, then for every object $A$ of ${\cal
A}$, $W_{A\underline{\hspace{1.3ex}}}$ is a lax functor and
$V_{A\underline{\hspace{1.3ex}}}$ is its rectification. As for
simplicial spaces, for a (lax) functor $W\!:\Delta^{op}\str\cat$,
we denote the unique arrow $W_0\str (W_1)^0$ by $p_0$, and when
$m\geq1$, we have $p_m=\langle i_1,\ldots,i_m\rangle\!:W_m\str
(W_1)^m$.

\vspace{2ex}

\df{4.2} We say that a lax functor $W\!:\Delta^{op}\str\cat$ is
\emph{Segal's}, when for every $m\geq 0$, $p_m\!: W_m\str(W_1)^m$
is the identity. We say that a lax functor
$W\!:(\Delta^{op})^n\str\cat$ is \emph{Segal's}, when for every
$l\in\{0,\ldots, n-1\}$ and every $k\geq 0$, the lax functor
$W_{\underbrace{\mbox{\scriptsize $1\ldots
1$}}_{l}\underline{\hspace{1.3ex}}k\ldots
k}\!:\Delta^{op}\str\cat$ is Segal's.

\vspace{2ex}

We denote by $B\!:\cat\str\top$ the \emph{classifying space}
functor, i.e., the composition $|\mbox{\hspace{1ex}}|\circ N$,
where $N\!:\cat\str \top^{\Delta^{op}}$ is the \emph{nerve}
functor.

\prop{Proposition 4.3}{If $W\!:\Delta^{op}\str\cat$ is Segal's lax
functor and $V$ is its rectification, then $B\circ V$ is Segal's
simplicial space.}

By Definitions 3.4 and 4.2, the following generalization of
Proposition~4.3 is easily obtained.

\prop{Corollary 4.4}{If $W\!:(\Delta^{op})^n\str\cat$ is Segal's
lax functor and $V$ is its rectification, then $B\circ V$ is
Segal's multisimplicial space.}

For Corollary~4.4, we conclude that $B\circ V$ is a PRmss as in
the proof of Proposition~4.3 given in the appendix.

\section{An application}

Let $\cal M$ be an $n$-fold strict monoidal category and let
$\WM\!:(\Delta^{op})^n\str\cat$ be the $n$-fold reduced bar
construction defined as in \cite{CP13}. The main result of that
paper says that $\WM$ is a lax functor and it is easy to verify
that it is Segal's. (Note that $\WM_{k_1\ldots k_l
\underline{\hspace{1.3ex}}\ldots k_{n-1}}$ is not Segal's when
$k_j>1$, for some $1\leq j\leq l$.)

For $V$ being a rectification of $\WM$, we have the following.

\prop{Theorem~5.1}{If $BV_{1\ldots 1}$, with respect to the
H-space structure obtained by Remark~3.5 when $l=n-1$, is
grouplike, then $B{\cal M}\simeq \Omega^n|B\circ V|$. }

\dkz By Corollary~4.4, we have that $B\circ V$ is Segal's
multisimplicial space. Hence, by Proposition~3.7, $BV_{1\ldots
1}\simeq \Omega^n|B\circ V|$. Since $V$ is a rectification of
$\WM$, by relying on Remark A1 of the appendix, we conclude that
$BV_{1\ldots 1}\simeq B\WM_{1\ldots 1}$. Together with the fact
that $\WM_{1\ldots 1}={\cal M}$, we obtain that
\[
B{\cal M}\simeq \Omega^n|B\circ V|.
\]

\vspace{-2.5em} \mbox{\hspace{2em}} \qed

This means that up to group completion (see \cite{S74} and
\cite{McDS76}), for every $n$-fold strict monoidal category $\cal
M$, the classifying space $B{\cal M}$ is an $n$-fold loop space.
When $\cal M$ contains a terminal or initial object, we have that
$B{\cal M}$, and hence $BV_{1\ldots 1}$, is path-connected. In
that case, by \cite[Proposition~8.4.4]{A11}, $BV_{1\ldots 1}$ is
grouplike, and $|B\circ V|$ is an $n$-fold delooping of~$B{\cal
M}$.

\section{Appendix}

\noindent{\sc Proof of Lemma 3.1.} First, we prove that $\langle
X_1,m,\ast\rangle$ is an H-space. Let $j_1\!:X_1\str X_1\times
X_1$ be such that $j_1(x)=(x,\ast)$, and analogously, let
$j_2\!:X_1\str X_1\times X_1$ be such that $j_2(x)=(\ast,x)$. By
the assumption, $X_0$ is contractible. Hence, $d^1_0$ is homotopic
to the constant map to $x_0$ and therefore $s^1_0\circ d^1_0$ is
homotopic to the constant map to $\ast$. We conclude that

\begin{tabbing}
\hspace{1.5em}$j_1\simeq \langle\mj_{X_1},s^1_0\circ d^1_0\rangle
= \langle d^2_2\circ s^2_1,d^2_0\circ s^2_1\rangle = \langle
d^2_2, d^2_0\rangle \circ s^2_1 = p_2\circ s^2_1$,
\end{tabbing}
i.e.,\ $p_2^{-1}\circ j_1\simeq s^2_1$. Hence,
\begin{tabbing}
\hspace{1.5em}$m\circ j_1= d^2_1\circ p_2^{-1}\circ j_1 \simeq
d^2_1\circ s^2_1 = \mj_{X_1}$.
\end{tabbing}
Analogously, $m\circ j_2\simeq \mj_{X_1}$ and we have that
$\langle X_1,m,\ast\rangle$ is an H-space.

\vspace{1ex}

Next, we prove that $m$ is associative up to homotopy, i.e., that
\[
m\circ(m\times\mj)\simeq m\circ(\mj\times m).
\]
Consider $p_3\!:X_3\str (X_1)^3$ for which we have:

\begin{tabbing}
\hspace{1.5em}\=$p_3$ \=$= \langle\langle
i_1,i_2\rangle,i_3\rangle = \langle\langle d^2_2\circ
d^3_3,d^2_0\circ d^3_3\rangle,i_3\rangle = \langle p_2\circ
d^3_3,i_3\rangle$
\\[1ex]
\>\>$=(p_2\times\mj)\circ\langle d^3_3,i_3\rangle$,\quad and
analogously
\\[1.5ex]
\> $p_3=(\mj\times p_2)\circ\langle i_1, d^3_0\rangle$.
\end{tabbing}

Since $p_2$ and $p_3$ are homotopy equivalences, we have that
$\langle d^3_3,i_3\rangle$ and $\langle i_1,d^3_0\rangle$ are
homotopy equivalences, too. Moreover,

\begin{tabbing}
\hspace{1.5em}\=$(1)$\quad \=$\langle d^3_3,i_3\rangle ^{-1}\simeq
p_3^{-1}\circ (p_2\times\mj)$, and
\\[1ex]
\> $(2)$ \> $\langle i_1,d^3_0\rangle^{-1}\simeq p_3^{-1}\circ
(\mj\times p_2)$.
\end{tabbing}

Also, we show that

\begin{tabbing}
\hspace{1.5em}\=$(3)$\quad \=$d^2_1\times\mj \simeq p_2\circ
d^3_1\circ p_3^{-1}\circ (p_2\times\mj)$, and
\\[1ex]
\> $(4)$ \> $\mj \times d^2_1 \simeq p_2\circ d^3_2\circ
p_3^{-1}\circ (\mj\times p_2)$.
\end{tabbing}

We have
\begin{tabbing}
\hspace{1.5em}\=$(d^2_1\times\mj)\circ \langle d^3_3, i_3\rangle$
\=$=\langle d^2_1\circ d^3_3,i_3 \rangle =\langle d^2_2\circ
d^3_1, d^2_0\circ d^3_1\rangle = \langle d^2_2, d^2_0\rangle \circ
d^3_1$
\\[1ex]
\> \> $= p_2\circ d^3_1$,
\end{tabbing}
which together with $(1)$ delivers $(3)$. Also,

\begin{tabbing}
\hspace{1.5em}\=$(\mj\times d^2_1)\circ \langle i_1,d^3_0\rangle$
\=$=\langle i_1, d^2_1\circ d^3_0 \rangle =\langle d^2_2\circ
d^3_2, d^2_0\circ d^3_2\rangle = \langle d^2_2, d^2_0\rangle \circ
d^3_2$
\\[1ex]
\> \> $= p_2\circ d^3_2$,
\end{tabbing}
which together with $(2)$ delivers $(4)$. Finally, we have
\begin{tabbing}
\hspace{1.5em}\=$m\circ(m\times\mj)$ \=$=d^2_1\circ p_2^{-1}\circ
(d^2_1\times \mj)\circ (p_2^{-1}) \simeq d^2_1\circ d^3_1\circ
p_3^{-1}$, by $(3)$
\\[1ex]
\> \> $= d^2_1\circ d^3_2\circ p_3^{-1} \simeq d^2_1\circ
p_2^{-1}\circ (\mj\times d^2_1)\circ (\mj\times p_2^{-1})$, by
$(4)$
\\[1ex]
\> \> $= m\circ(\mj\times m)$. \` $\dashv$
\end{tabbing}

\vspace{3ex}

\noindent{\sc Proof of Lemma 3.3.} Let $Z\!:\Delta^{op}\str\top$
be $Y^{(1)}$. By Remark~2.7 (a), it is a PRmss. We have to show
that for every $m\geq 0$, the map $p_m\!:Z_m\str (Z_1)^m$ is a
homotopy equivalence.

Let $m=0$ and let $T$ be the trivial simplicial space with
$T_k=\{\ast\}$. Consider the simplicial space
$Y_{0\underline{\hspace{1.3ex}}}\!:\Delta^{op}\str\top$, which is
a PRmss by Remark 2.7~(c). By Proposition~2.9 and Lemma~2.10, both
$T$ and $Y_{0\underline{\hspace{1.3ex}}}$ are proper. The
following simplicial map is obtained by the assumptions (the
diagrams are commutative since $\{\ast\}$ is terminal).
\begin{center}
\begin{picture}(210,60)(20,55)

\put(50,85){\makebox(0,0){$\downarrow$}}
\put(55,87){\makebox(0,0)[l]{$\simeq$}}

\put(130,85){\makebox(0,0){$\downarrow$}}
\put(135,87){\makebox(0,0)[l]{$\simeq$}}

\put(210,85){\makebox(0,0){$\downarrow$}}
\put(215,87){\makebox(0,0)[l]{$\simeq$}}

\put(20,105){\makebox(0,0){$\ldots$}}

\put(50,105){\makebox(0,0){$Y_{02}$}}

\put(90,115){\makebox(0,0){$\rightarrow$}}
\put(90,105){\makebox(0,0){$\rightarrow$}}
\put(90,95){\makebox(0,0){$\rightarrow$}}

\put(90,110){\makebox(0,0){$\leftarrow$}}
\put(90,100){\makebox(0,0){$\leftarrow$}}

\put(130,105){\makebox(0,0){$Y_{01}$}}

\put(170,110){\makebox(0,0){$\rightarrow$}}
\put(170,100){\makebox(0,0){$\rightarrow$}}

\put(170,105){\makebox(0,0){$\leftarrow$}}

\put(210,105){\makebox(0,0){$Y_{00}$}}

\put(20,65){\makebox(0,0){$\ldots$}}

\put(50,65){\makebox(0,0){$\{\ast\}$}}

\put(90,67){\makebox(0,0){$\rightarrow$}}

\put(90,62){\makebox(0,0){$\leftarrow$}}

\put(130,65){\makebox(0,0){$\{\ast\}$}}

\put(170,67){\makebox(0,0){$\rightarrow$}}

\put(170,62){\makebox(0,0){$\leftarrow$}}

\put(210,65){\makebox(0,0){$\{\ast\}$}}

\put(-20,105){\makebox(0,0){$Y_{0\underline{\hspace{1.3ex}}}:$}}
\put(-20,65){\makebox(0,0){$T:$}}

\end{picture}
\end{center}
By Lemma~2.11, we have that
$|Y_{0\underline{\hspace{1.3ex}}}|\simeq |T|=\{\ast\}$ via the
unique map. Since $Z_0=|Y_{0\underline{\hspace{1.3ex}}}|$ and
$(Z_1)^0=\{\ast\}$, we are done.

\vspace{1ex}

Let $m>0$. Consider the simplicial spaces
$Y_{m\underline{\hspace{1.3ex}}}$ and
$(Y_{1\underline{\hspace{1.3ex}}})^m$, which are proper by Remark
2.7~(c), Proposition~2.9, Lemma~2.10 and Remark~2.13. The
following simplicial map is obtained by the assumptions (it is
straightforward to verify that the diagrams are commutative).
\begin{center}
\begin{picture}(210,60)(20,55)

\put(50,85){\makebox(0,0){$\downarrow$}}
\put(55,87){\makebox(0,0)[l]{$\simeq$}}

\put(130,85){\makebox(0,0){$\downarrow$}}
\put(135,87){\makebox(0,0)[l]{$\simeq$}}

\put(210,85){\makebox(0,0){$\downarrow$}}
\put(215,87){\makebox(0,0)[l]{$\simeq$}}

\put(20,105){\makebox(0,0){$\ldots$}}

\put(50,105){\makebox(0,0){$Y_{m2}$}}

\put(90,115){\makebox(0,0){$\rightarrow$}}
\put(90,105){\makebox(0,0){$\rightarrow$}}
\put(90,95){\makebox(0,0){$\rightarrow$}}

\put(90,110){\makebox(0,0){$\leftarrow$}}
\put(90,100){\makebox(0,0){$\leftarrow$}}

\put(130,105){\makebox(0,0){$Y_{m1}$}}

\put(170,110){\makebox(0,0){$\rightarrow$}}
\put(170,100){\makebox(0,0){$\rightarrow$}}

\put(170,105){\makebox(0,0){$\leftarrow$}}

\put(210,105){\makebox(0,0){$Y_{m0}$}}

\put(20,65){\makebox(0,0){$\ldots$}}

\put(50,65){\makebox(0,0){$(Y_{12})^m$}}

\put(90,75){\makebox(0,0){$\rightarrow$}}
\put(90,65){\makebox(0,0){$\rightarrow$}}
\put(90,55){\makebox(0,0){$\rightarrow$}}

\put(90,70){\makebox(0,0){$\leftarrow$}}
\put(90,60){\makebox(0,0){$\leftarrow$}}

\put(130,65){\makebox(0,0){$(Y_{11})^m$}}

\put(170,70){\makebox(0,0){$\rightarrow$}}
\put(170,60){\makebox(0,0){$\rightarrow$}}

\put(170,65){\makebox(0,0){$\leftarrow$}}

\put(210,65){\makebox(0,0){$(Y_{10})^m$}}

\put(-20,105){\makebox(0,0){$Y_{m\underline{\hspace{1.3ex}}}:$}}
\put(-20,65){\makebox(0,0){$(Y_{1\underline{\hspace{1.3ex}}})^m:$}}

\end{picture}
\end{center}
By Lemma~2.11, we have that
\[
|\langle Y(i_1,\underline{\hspace{1.3ex}}),\ldots,Y(i_m,
\underline{\hspace{1.3ex}})\rangle|\!:|Y_{m\underline{\hspace{1.3ex}}}|
\str |(Y_{1\underline{\hspace{1.3ex}}})^m|
\]
is a homotopy equivalence. Also, for \top, the realization functor
$|\mbox{\hspace{1ex}}|$ preserves products (see
\cite[Theorem~14.3]{M67}, \cite[III.3, Theorem]{GZ67} and
\cite[Corollary~11.6]{M72}). Namely, for
$\pi_k\!:(Y_{1\underline{\hspace{1.3ex}}})^m \str
Y_{1\underline{\hspace{1.3ex}}}$, $1\leq k\leq m$ being the $k$th
projection,
\[
\langle|\pi_1|,\ldots,|\pi_m|\rangle\!:
|(Y_{1\underline{\hspace{1.3ex}}})^m| \str
|Y_{1\underline{\hspace{1.3ex}}}|^m
\]
is a homeomorphism (${|\mbox{\hspace{1ex}}|}$ is strong monoidal;
see \cite[Example~6.2.2]{R14}). Hence,
\[
\langle|\pi_1|,\ldots,|\pi_m|\rangle\circ |\langle
Y(i_1,\underline{\hspace{1.3ex}}),\ldots,Y(i_m,\underline{\hspace{1.3ex}})\rangle|\!:
|Y_{m\underline{\hspace{1.3ex}}}| \str
|Y_{1\underline{\hspace{1.3ex}}}|^m
\]
is a homotopy equivalence.

The following easy computation, in which $\langle
Y(i_1,\underline{\hspace{1.3ex}}),\ldots,Y(i_m,\underline{\hspace{1.3ex}})\rangle$
is abbreviated by $\alpha$,

\begin{tabbing}
\hspace{1.5em}\=$\langle|\pi_1|,\ldots,|\pi_m|\rangle \circ
|\alpha|$ \=$=\langle|\pi_1|\circ |\alpha|,\ldots,|\pi_m|\circ
|\alpha|\rangle = \langle|\pi_1\circ \alpha|,\ldots,|\pi_m\circ
\alpha|\rangle$
\\[1ex]
\> \> $=
\langle|Y(i_1,\underline{\hspace{1.3ex}})|,\ldots,|Y(i_m,\underline{\hspace{1.3ex}})|\rangle
= \langle Z(i_1),\ldots,Z(i_m)\rangle$
\end{tabbing}
shows that the map $p_m=\langle Z(i_1),\ldots,Z(i_m)\rangle$ is a
homotopy equivalence between
$Z_m=|Y_{m\underline{\hspace{1.3ex}}}|$, and
$(Z_1)^m=|Y_{1\underline{\hspace{1.3ex}}}|^m$. \qed

\vspace{3ex}

\noindent{\sc Some preliminary remarks for Proposition 4.3.} Let
\textbf{2} be the category with two objects (0 and 1) and one
nonidentity arrow $h\!:0\str 1$. Let $I_0,I_1\!:{\cal C}\str {\cal
C}\times \textbf{2}$ be the functors such that for every object
$C$ of $\cal C$, we have that $I_0(C)=(C,0)$ and $I_1(C)=(C,1)$.
Let $F,G\!:{\cal C}\str{\cal D}$ be two functors. There is a
bijection between the set of natural transformations
$\alpha\!:F\strt G$, and the set of functors $A\!:{\cal C}\times
\textbf{2}\str {\cal D}$ such that $A\circ I_0=F$ and $A\circ
I_1=G$. This bijection maps $\alpha\!:F\strt G$ to  $A\!:{\cal
C}\times \textbf{2}\str {\cal D}$ such that
\[
A(C,0)=FC,\quad A(C,1)=GC,\quad A(f,\mj_0)=Ff,\quad A(f,\mj_1)=Gf,
\]
and for $f\!:C\str C'$,
\[
A(f,h)=Gf\circ\alpha_C=\alpha_{C'}\circ Ff.
\]
Its inverse maps $A\!:{\cal C}\times \textbf{2}\str {\cal D}$ to
$\alpha\!:F\strt G$ such that $\alpha_C=A(\mj_C,h)$.

The nerve functor $N$ preserves products on the nose, hence, the
classifying space functor $B=|\mbox{\hspace{1ex}}|\circ N$
preserves products too. Therefore, the spaces $B{\cal C}\times I$
(i.e., $B{\cal C}\times B\textbf{2}$) and $B({\cal C}\times
\textbf{2})$ are homeomorphic and we have the following.

\prop{Remark A1}{Every natural transformation $\alpha\!:F\strt G$
gives rise to the homotopy
\[
B{\cal C}\times I\stackrel{\approx}{\longrightarrow} B({\cal
C}\times \emph{\textbf{2}}) \stackrel{BA\;}{\longrightarrow}
B{\cal D}
\]
between the maps $BF$ and $BG$.}


\noindent{\sc Proof of Proposition 4.3.} By the isomorphism
mentioned in Definition~2.2, we have that $N\circ V$ corresponds
to a multisimplicial set $X\!:\Delta^{op}\times
\Delta^{op}\str\top$ and $B\circ V$ is $X^{(1)}$. Hence, it is a
PRmss.

We have to show that for every $m\geq 0$, $p_m\!:BV_m\str
(BV_1)^m$ is a homotopy equivalence, where we denote again by
$p_0$ the unique map from $BV_0$ to $(BV_1)^0$ and by $p_m$ the
map $\langle BV(i_1),\ldots, BV(i_m) \rangle$.

When $m=0$, we show that $BJ_0\!:BW_0\str BV_0$ is a homotopy
inverse to $p_0$. Since $W_0$ and $(V_1)^0$ are the same trivial
category and $BW_0=(BV_1)^0=\{\ast\}$, it is easy to conclude that
$p_0\circ BJ_0\simeq \mj_{(BV_1)^0}$, and that $p_0=BE_0$. The
latter, by the adjunction $J_0\dashv E_0$ and Remark A1, delivers
$BJ_0\circ p_0\simeq \mj_{BV_0}$.

When $m\geq 1$, we have for every $1\leq j\leq m$, the following
natural transformations.

\begin{center}
\begin{picture}(320,80)
\put(0,40){\makebox(0,0){$V_m$}} \put(60,10){\makebox(0,0){$V_1$}}
\put(120,40){\makebox(0,0){$W_1$}}
\put(60,70){\makebox(0,0){$W_m$}}
\put(25,20){\makebox(0,0)[t]{$i_j$}}
\put(25,60){\makebox(0,0)[b]{$E_m$}}
\put(95,60){\makebox(0,0)[b]{$i_j$}}
\put(95,20){\makebox(0,0)[t]{$E_1$}}

\put(60,40){\makebox(0,0){$\Downarrow$}}
\put(65,40){\makebox(0,0)[l]{\scriptsize $E_{i_j}$}}

\put(10,35){\vector(2,-1){40}} \put(10,45){\vector(2,1){40}}
\put(70,65){\vector(2,-1){40}} \put(70,15){\vector(2,1){40}}

\put(200,40){\makebox(0,0){$W_m$}}
\put(260,10){\makebox(0,0){$W_1$}}
\put(320,40){\makebox(0,0){$V_1$}}
\put(260,70){\makebox(0,0){$V_m$}}
\put(225,20){\makebox(0,0)[t]{$i_j$}}
\put(225,60){\makebox(0,0)[b]{$J_m$}}
\put(295,60){\makebox(0,0)[b]{$i_j$}}
\put(295,20){\makebox(0,0)[t]{$J_1$}}

\put(260,40){\makebox(0,0){$\Uparrow$}}
\put(265,37){\makebox(0,0)[l]{\scriptsize $J_{i_j}$}}

\put(210,35){\vector(2,-1){40}} \put(210,45){\vector(2,1){40}}
\put(270,65){\vector(2,-1){40}} \put(270,15){\vector(2,1){40}}

\end{picture}
\end{center}
By using the monoidal structure of \cat\ given by 2-products and
the fact that $\langle i_1,\ldots,i_m\rangle\!:W_m\str (W_1)^m$ is
the identity, we obtain the following two natural transformations.

\begin{center}
\begin{picture}(320,80)
\put(0,70){\makebox(0,0){$V_m$}}
\put(60,10){\makebox(0,0){$(V_1)^m$}}
\put(120,70){\makebox(0,0){$W_m$}}
\put(5,40){\makebox(0,0)[t]{$\langle i_1,\ldots,i_m\rangle$}}
\put(60,75){\makebox(0,0)[b]{$E_m$}}
\put(105,40){\makebox(0,0)[t]{$(E_1)^m$}}

\put(60,40){\makebox(0,0){$\Downarrow$}}
\put(60,50){\makebox(0,0){\scriptsize $\langle E_{i_1},\ldots,
E_{i_m}\rangle$}}

\put(10,60){\vector(1,-1){40}} \put(10,70){\vector(1,0){97}}
\put(70,20){\vector(1,1){40}}

\put(200,10){\makebox(0,0){$W_m$}}
\put(320,10){\makebox(0,0){$(V_1)^m$}}
\put(260,70){\makebox(0,0){$V_m$}}
\put(285,45){\makebox(0,0)[bl]{$\langle i_1,\ldots,i_m\rangle$}}
\put(260,5){\makebox(0,0)[t]{$(J_1)^m$}}
\put(230,45){\makebox(0,0)[br]{$J_m$}}

\put(260,40){\makebox(0,0){$\Uparrow$}}
\put(260,30){\makebox(0,0){\scriptsize $\langle J_{i_1},\ldots,
J_{i_m}\rangle$}}

\put(270,60){\vector(1,-1){40}} \put(210,10){\vector(1,0){97}}
\put(210,20){\vector(1,1){40}}

\end{picture}
\end{center}

For $\pi_k\!:{\cal C}^m \str {\cal C}$, $1\leq k\leq m$ being the
$k$th projection,
\[
\langle B\pi_1,\ldots,B\pi_m\rangle\!:B{\cal C}^m\str (B{\cal
C})^m
\]
is a homeomorphism whose inverse we denote by $q_m({\cal C})$. It
is easy to verify that for $F,F_1,\ldots,F_m\!:{\cal C}\str{\cal
D}$ we have $B\langle F_1,\ldots,F_m\rangle=q_m({\cal D})\langle
BF_1,\ldots,BF_m\rangle$ and $BF^m\circ q_m({\cal C})=q_m({\cal
D})\circ (BF)^m$.

By Remark~A1, the transformations mentioned above give rise to
\begin{tabbing}
\hspace{1.5em}\=$(\dagger)$\hspace{3em} \=$BE_m$\hspace{1em} \=
$\simeq B(E_1)^m\circ B\langle V(i_1),\ldots,V(i_m)\rangle$
\\[1ex]
\>\>\> $= q_m(W_1)\circ (BE_1)^m\circ \langle
BV(i_1),\ldots,BV(i_m)\rangle$
\\[1ex]
\>\>\> $= q_m(W_1)\circ (BE_1)^m\circ p_m$, and
\\[2ex]
\>$(\dagger\dagger)$\>$B(J_1)^m$\> $\simeq B\langle
V(i_1),\ldots,V(i_m)\rangle \circ BJ_m$
\\[1ex]
\>\>\> $= q_m(V_1)\circ \langle BV(i_1),\ldots,BV(i_m)\rangle
\circ BJ_m$
\\[1ex]
\>\>\> $= q_m(V_1)\circ p_m \circ BJ_m$.
\end{tabbing}

The following calculation shows that
\[
BJ_m\circ q_m(W_1)\circ (BE_1)^m\!:(BV_1)^m\str BV_m
\]
is a homotopy inverse to $p_m$.
\begin{tabbing}
\hspace{1.5em}\=$\mj_{BV_m}$\hspace{1.2em}\= $\simeq BJ_m\circ
BE_m$,\quad by $J_m\dashv E_m$, Remark~A1
\\[1ex]
\>\> $\simeq BJ_m\circ q_m(W_1)\circ (BE_1)^m\circ p_m$,\quad by
$(\dagger)$
\\[2ex]
\>$\mj_{(BV_1)^m}$\> $\simeq q_m^{-1}(V_1)\circ B(J_1)^m\circ
B(E_1)^m \circ q_m(V_1)$,\quad by $J_1\dashv E_1$, Remark~A1
\\[1ex]
\>\> $\simeq p_m\circ BJ_m\circ q_m(W_1)\circ (BE_1)^m$,\quad by
$(\dagger\dagger)$. \`$\dashv$
\end{tabbing}

\pagebreak

\end{document}